 \documentclass[authoryear,3p, 9pt,times,final,twocolumn]{elsarticle} 

\usepackage{amsfonts}
\usepackage[htt]{hyphenat}

\usepackage{eulervm}
\usepackage[mathscr, mathcal]{euscript}
\usepackage{anyfontsize}
\usepackage{float}
\restylefloat{figure}
\restylefloat{table}
\usepackage{graphicx}
\usepackage{caption}
\usepackage{subcaption}
\usepackage{amsmath}
\allowdisplaybreaks

\usepackage{amsmath}

\usepackage[normalem]{ulem}

\usepackage{float}
 \usepackage{graphics}
\usepackage[bookmarks=false]{hyperref}

\hypersetup{%
  colorlinks=true,
  linkcolor=blue,
  pdfstartview=FitH,
  breaklinks=true, 
  breaklinks=true,
  bookmarksopen=false,
    filecolor=magenta,
    urlcolor=cyan,
    pdftitle={Overleaf Example},
    pdfpagemode=FullScreen,
  unicode   =false,    
  pdftoolbar=false,    
  pdfmenubar=true      
  }
\usepackage{bookmark}

\usepackage{longtable,lscape,verbatim}
\usepackage{amsmath}
\usepackage{graphicx}
\usepackage[font=small,labelfont=bf]{caption}
\usepackage{latexsym}
\usepackage{amsmath}
\usepackage{todonotes}
\usepackage[nomargin,inline,marginclue,draft]{fixme}
\makeatletter
\renewcommand*\FXLayoutInline[3]{%
	{\@fxuseface{inline}\ignorespaces[#3 \fxnotename{#1}: #2]}}
\makeatother
\FXRegisterAuthor{sv}{asv}{\color{red}Me}
\allowdisplaybreaks
\graphicspath{{figs/}}
\usepackage{verbatim}
\usepackage{amsmath,amsthm}
\allowdisplaybreaks
\usepackage{algorithm}
\usepackage{algpseudocode}
\usepackage{fancybox}
 \usepackage{times}
\usepackage{fixme}

\newtheorem{definition}{Definition}

\newtheorem{thm}{Theorem}
\newtheorem{lem}[thm]{Lemma}

\usepackage{amsmath}
\usepackage{etoolbox}
\AtBeginEnvironment{align}{\setcounter{subeqn}{0}}
\newcounter{subeqn} %

\usepackage{listings}
\usepackage{charter}
\usepackage[expert]{mathdesign}
\usepackage{pgfplots}
\usepackage{times} 
\usepackage[compatibility=false]{caption}

\usepackage{graphicx}
\usepackage{caption}
\usepackage{subcaption}

\usepackage{filecontents}
\usepackage{pgf, tikz}
\usetikzlibrary{arrows, automata}

\usepackage{lipsum}

\allowdisplaybreaks[1] 

\tikzstyle{vertex}=[circle,fill=black!25,minimum size=20pt,inner sep=0pt]
\tikzstyle{selected vertex} = [vertex, fill=red!24]
\tikzstyle{blue vertex} = [vertex, fill=blue!24]
\tikzstyle{edge} = [draw,thick,-]
\tikzstyle{weight} = [font=\small]
\tikzstyle{selected edge} = [draw,line width=5pt,-,red!50]
\tikzstyle{ignored edge} = [draw,line width=5pt,-,black!20]

\usepackage[pdf]{pstricks}%
\usepackage{pstricks-add, pst-grad}
\definecolor[ps]{bordercolor}{rgb}{0.06 0.42 0.11}
\definecolor[ps]{greenbegin}{rgb}{0.47 0.94  0.151}
\definecolor[ps]{greenend}{rgb}{0.23 0.65  0.02}

\usepackage{setspace}
\allowdisplaybreaks

\usepackage{todonotes}

\tikzset{inlinenotestyle/.append style={align=justify}}

\usepackage{lipsum} 

\usepackage{color}

\usepackage{xcolor}

\usepackage{dsfont}
\usepackage{longtable}
\usepackage{lscape}
\usepackage{ctable}

\DeclareMathAlphabet{\pazocal}{OMS}{zplm}{m}{n}

\usepackage[utf8]{inputenc}

\usepackage{cancel}
\makeatletter
\newcommand{\crossout}[1]{%
  \begingroup
  \settowidth{\dimen@}{#1}%
  \setlength{\unitlength}{0.05\dimen@}%
  \settoheight{\dimen@}{#1}%
  \count@=\dimen@
  \divide\count@ by \unitlength
  \begin{picture}(0,0)
  \put(0,0){\line(20,\count@){20}}
  \put(0,\count@){\line(20,-\count@){20}}
  \end{picture}%
  #1%
  \endgroup
}
\makeatother

\usepackage{lscape}
\DeclareGraphicsExtensions{.png,.eps}

\usepackage{pdflscape}

\usepackage{tikz-network}
\usetikzlibrary{shapes}

\usepackage{xcolor,cancel}

\usepackage{ctable} 

\usepackage{mathdesign}
\usepackage[OMLmathrm,OMLmathbf]{isomath}

\journal{}

\begin{document}
\begin{frontmatter}



\title{A Neural Benders Decomposition for the Hub Location Routing Problem}
\author[CRISTAL,LILLE,ITBA,SPG]{Rahimeh Neamatian Monemi}
\author[LGI2A,cor2]{Shahin Gelareh}
%
%
\corref{cor1}
\cortext[cor2]{Corresponding author, Shahin Gelareh, shahin.gelareh@univ-artois.fr}

\address[LILLE]{Université de Lille, France}
\address[ITBA]{IT and Business Analytics Ltd, UK}
\address[SPG]{Sharkey Predictim Globe, France}
\address[LGI2A]{Département R\&T, IUT de Béthune, Université d'Artois, F-62000 Béthune, France}

\begin{abstract}
In this study, we propose an imitation learning framework designed to enhance the Benders decomposition method. Our primary focus is addressing degeneracy in subproblems with multiple dual optima, among which Magnanti-Wong technique identifies the non-dominant solution. We develop two policies. In the first policy, we replicate the Magnanti-Wong method and learn from each iteration. In the second policy, our objective is to determine a trajectory that expedites the attainment of the final subproblem dual solution. We train and assess these two policies through extensive computational experiments on a network design problem with flow subproblem, confirming that the presence of such learned policies significantly enhances the efficiency of the decomposition process.
\end{abstract}


\end{frontmatter}
%
\normalsize{ }
\section{Introduction}
Benders decomposition \cite{benders1962partitioning,benders2005partitioning} serves as a method for tackling mixed integer linear programming (MILP) problems, proving effective for both deterministic and stochastic instances of mixed integer programming. The fundamental concept behind Benders decomposition is to project the problem onto the domain of complicating variables, leading to a relaxation of the original problem. This outer approximation is then refined by introducing cuts derived from the extreme points of a linear programming polytope. These cuts are parameterized within the feasibility domain of the complicating variable concerning the initial problem. Essentially, Benders decomposition divides a given MILP problem into two components: a Master Problem (MP) featuring the integer variables (also known as complicating variables), and a Subproblem (SP) parameterized by the variables from the master problem. In its simplest form, the subproblem in Benders decomposition takes the shape of a linear programming problem.  Although in some cases such as in some stochastic programming models and of course in combinatorial Benders \cite{codato2006combinatorial}, one may encounter  (mixed) integer programming subproblems.\\

The Master Problem (MP), which constitutes a relaxation, typically presents a smaller problem for which we possess an efficient algorithm, or it is amenable to solution using off-the-shelf solvers, proving considerably more efficient compared to addressing the entire problem. A similar rationale extends to the subproblem, which could even be partitioned into smaller subproblems under certain circumstances. Assuming a feasible model is available, the iterative procedure proceeds to refine this outer approximation in the following manner: The MP proposes a (potentially optimal) solution, and the subproblem faces one of the two scenarios: 1) The subproblem is feasible given the parameters provided by the master problem. In this case, an extreme point of the subproblem polytope emerges as an \emph{optimality cut}, effectively truncating the current master solution and enhancing the outer approximation. 2) The subproblem becomes infeasible when employing the parameters conveyed by the master problem. In such instances, the subproblem returns a \emph{feasibility cut} along with a certificate of infeasibility. This prohibits the master problem from reiterating the same information and encourages it to supply alternative insights. Each iteration yields a lower bound resulting from the master problem---functioning as a relaxation of the original problem---and an upper bound achieved through a combination of the master problem and the subproblem solutions. Termination of the iterations occurs when the optimality cuts are no longer viable for the current master problem or when the disparity between the lower and upper bounds falls below a specified threshold.\\

Given the following  mixed integer programming:

\begin{align}
z^* = \min_{x,y} \quad & c^T x + d^T y \label{p:obj}\\
\text{s.t.}~~ u: ~~\quad & Ax + By \leq b \label{p:eq1}\\
           ~~ v:~~ \quad & Fx \leq g \label{p:eq2}\\
& \mathbf{x} \in \mathbb{Z}, \mathbf{y} \in \mathbb{R}\label{p:eq3}
\end{align}
where $x$ and $y$ are the decision variables, $c$ and $d$  are coefficient vectors, $A$, $B$ and $F$ are coefficient matrices of appropriate size, and $b$ and $g$ are constraint vectors. Furthermore, we define $\mathcal{X} = \{x: Fx\leq g, x\in \mathbb{Z} \}$.\\

The Master Problem (MP) follows:

\begin{align}
z^*_{MP} = \min_{x} \quad & c^T x + \eta \label{mp:obj}\\
\text{s.t.}  \quad & Fx \leq g \label{mp:eq1}\\
             \quad &  \eta \geq \overline{u} (b- A{x})  & \text{\small for all extreme points of SP} \label{mp:eq2}\\
             \quad &  0 \geq \overline{u} (b- A{x})   & \text{\small for all extreme rays of SP} \label{mp:eq3}\\
 \quad & \mathbf{x}  \in \mathbb{Z}, \eta \geq 0\label{mp:eq4}
\end{align}

where $\mathbf{u}$ is the vector of duals associated to the constraints \eqref{p:eq1} and $\eta$ is the variable linking the master problem to the following Subproblem (SP):\

\begin{align}
z^*_{SP(y)} = \min_{y} \quad &d^T y \label{sp:obj}\\
\text{s.t.}~~ u: ~~\quad & B~y \leq b - A\overline{x}  \label{sp:eq1}\\
&  {y} \in \mathbb{R}\label{sp:eq2}
\end{align}
with the following Subproblem  Dual (SPD):
\begin{align}
z^*_{SPD(u)} = \max \quad & u(b- A\overline{x}) \label{spd:obj}\\
\text{s.t.} \quad & u~B \geq d  \label{spd:eq1}\\
&   {u}\in \mathbb{R}^+\cup\{0\} \label{spd:eq2}
\end{align}

The iterative process of the textbook Benders decomposition is illustrated in  \autoref{alg:Benders}. An efficient and modern implementation of Benders decomposition involves separating and adding the Benders cuts during the branch-and-bound process used to solve the master problem. This approach proves particularly advantageous for mid-sized to large-scale Mixed-Integer Linear Programming (MILP) problems, as it significantly reduces computational overhead. Unlike the requirement to solve the master problem to optimality and generate a cut for the optimal solution every time, this approach allows for generating Benders cuts for both integer and fractional solutions encountered at various nodes along the branch-and-bound tree of the master problem. Moreover, when usinh SPD instead of SP, further computational efficiency can be achieved. This is due to the SPD polytope's independence from the variable values  informed by the master problem. Therefore, any new information from the master problem only updates the SPD's objective function. This enables a warm-start from the previous optimal basis, and a few pivoting in simplex suffices to obtain the new optimal solution, as feasibility is already established.\\

\begin{algorithm}
\small
  \caption{Benders Decomposition for Deterministic Mixed Integer Programming}
  \begin{algorithmic}[1]
    \State Initialize the iteration counter $k=0$.
    \State Solve the  master problem (MP) using any available optimization solver, $z^k_{MP}$.
    \State $LB = -\infty$, $UB = \infty$ and $\epsilon = 0.001$
    \Repeat
      \State Solve the subproblem (SP) using the current MP solution $x^k$ to obtain the optimal values $u^k$.
      \If {SP is Feasible}
          \State Compute the optimality cut: $\eta \geq \overline{u} (b- A{x}) $.
      \Else
          \State Compute the feasibility cut: $0 \geq \overline{u} (b- A{x}) $.
      \EndIf
      \State Increment the iteration counter: $k \gets k+1$.
      \State Add the cuts to the master problem and re-solve it, $z^k_{MP}$.
      \State $LB = z^k_{MP}, UB =  z^k_{MP} + z^k_{SPD} - \eta$
      \State $gap = \frac{UB-LB}{UB}\times 100$
    \Until{$gap<\epsilon$}
  \end{algorithmic}\label{alg:Benders}
\end{algorithm}

This work assumes a deterministic problem and focuses on generating cuts from a single subproblem. However, it can readily accommodate scenarios involving multiple subproblems. We make the further assumption that the subproblem remains feasible for any solution derived from the master problem (referred to as complete recourse in stochastic programming). Consequently, we are able to omit the inequalities \eqref{mp:eq3} corresponding to the Benders feasibility cuts.\\ 

An influential factor affecting the efficiency of Benders decomposition, regardless of whether it is the traditional or modern implementation, is the cut selection strategy. This challenge becomes more pronounced when the subproblem exhibits degeneracy, leading to multiple optimal solutions within its dual. As a result, certain dual solutions may yield cuts that are dominated by others. An incorrect choice in this context would necessitate more Benders iterations to meet termination criteria, resulting in prolonged convergence. Thus, careful consideration is imperative in selecting a dual solution that generates precise cuts, contributing to a reduction in the iteration count. 

\section{Neural Benders Decomposition}\label{sec::solution}

The problem can be reformulated as a Markov Decision Process (MDP) and learn a policy for its solution. Transforming a decision problem into a Markov Decision Process (MDP) involves defining several key components. Firstly, at each iteration $t$, the state denoted as $s_t$ encapsulates all relevant information at the current node of the branch-and-bound tree. This information encompasses the solution pool, all previously visited incumbents, and integer solutions identified by heuristics within the Master Problem (MP). The terminal or final state is attained when either the effectiveness of Benders cuts in separating the integer node diminishes, or when Benders iterations meet termination criteria, often indicated by close proximity of the lower and upper bounds. At time $t$, the action $a_t$ represents the subsequent cut to be introduced into the master problem. The transition probability $p(s' | s, a)$ is computed by applying the action $a$ (i.e., adding the cut) to the state $s$, subsequently solving the Master Problem.\\


In this study, our objective is to learn and mimic an \emph{expert} proficient in isolating the subsequent cut to trim the Master Problem (MP) polytope. This endeavor is undertaken with the dual purpose of achieving a swift convergence in terms of both time and the frequency of invoking the separation algorithm.
Before we proceed any further, it is necessary to introduce some notations.

\begin{definition}(Core point, \cite{magnanti1981accelerating})
Any point $x^0$  in the relative interior of the closure, $x^c$, is  referred to as a \emph{core point}.
 \end{definition}

\begin{definition}(Dominance, \cite[Definition 1]{papadakos2008practical} and \cite{magnanti1981accelerating})
Let $u_1$ and $u_2$ denote the optimal solutions of SP at a given iteration $t$. The Benders cut corresponding to the dual solution $u^1$ dominates the one associated with dual solution $u^2$, if $(b-y)^Tu_1 \geq (b-y)^Tu_2$, with strict inequality for at least one point $y$ from SPD. As the cuts stem  from  $u_1$ and $u_2$, we can conclude that $u_1$ dominates $u_2$.
 \end{definition}

\begin{definition}(Pareto-optimal, \cite{magnanti1981accelerating})
A cut is considered Pareto-optimal when it is not dominated by any other cut. Here, $u$ is termed  Pareto-optimal.
\end{definition}

\subsection{Expert Policy}
We proceed to learn and imitate two policies as outlined below:

\subsubsection{Magnanti-Wong}\label{policy:MW}

 Our expert relies on the Magnanti-Wong technique, concentrating particularly on subproblems that display degeneracy and present multiple optimal solutions. In such cases, our goal is to identify the most incisive one among them. With each proposition originating from the master problem, this approach utilizes a \emph{core point}, specifically relative to the current iteration's approximation of the Master Problem (MP). This differs from Papadokos's approach that uses an \emph{approximate core point}. The expert then proceeds to select a point along the optimality face of SPD, aiming to generate the sharpest possible Benders cut.\\

 The constraint \eqref{spd:relint:face} confines our search space to the optimality face. To generate non-dominated cuts or Pareto-optimal cuts, one can employ the following formulation:

\begin{align}
\max  \quad & u(b- A{x}^{core}_t) \label{spd:relint}\\
\text{s.t.} \quad & u~B \geq d  \label{spd:eq1}\\
& z^*_{SPD(u)} =  u(b- A\overline{x}) \label{spd:relint:face}\\
&   {u}\in \mathbb{R}^+\cup\{0\} \label{spd:eq2}.
\end{align}

Here, $u(b - A{x}^{core}_t)$ represents the Euclidean distance between the core point of the iteration ${x}^{core}_t$ and the objective function of SPD.


\begin{definition}(Pareto-optimal, \cite[Theorem 1]{magnanti1981accelerating} and \cite[Definition 4]{papadakos2008practical})
The optimal solution of \eqref{spd:relint}-\eqref{spd:eq2}, where the point $x_t^{core}$ serves as a core point, results in a Pareto-optimal cut.
\end{definition}

One of the challenges faced by this algorithm is the selection of an appropriate \emph{core point} and the iterative updating of it. The sought-after \emph{core point} must lie within the integer hull of the master problem approximation at iteration $t$. It is important to note that a point situated in the gap between the integer hull and the linear hull in the master problem at iteration $t$ does not necessarily result in a non-dominated cut. The following lemma provides insight into the concept of a core point:

\begin{lem}\label{lem:relint}
Let $X^t_{pool} = \{x^t_0, x_1^t, x_2^t, \dots, x^t_p\}$ represent the set of integer feasible solutions of the master problem at iteration $t$. The point $x^o = \frac{\sum_{i=i}^{p}x_p^t}{|X^t_{pool}|}$ qualifies as a core point.
\end{lem}

In \cite{papadakos2008practical} and \cite{hosseini2021deepest}, it is demonstrated that one need not limit movement to the optimality face of the SPD. Instead, it suffices to locate an optimal point that maximizes the distances \eqref{spd:relint}. Therefore, the constraint \eqref{spd:relint:face} can be omitted from the model \eqref{spd:relint}-\eqref{spd:eq2}. However, for the scope of our current work, it remains necessary to retain this constraint and employ the method presented in \cite{magnanti1981accelerating}. In particular, in \cite{hosseini2021deepest}, the authors extended the concept of distance to $l_p$ and demonstrated that various options exist for selecting a dual solution that appears to deliver the deepest cut, even if it is not an optimal dual. Consequently, this approach does not deal with the degeneracy not it depend on the notion of the optimality face or the optimal solution of dual.


\subsubsection{Last Iteration Duals}\label{policy:lid}

Based our assumption, an optimal  $u^*$ always exists, a point where optimality has been conclusively established. In the realm of linear programming and column generation, \cite{babaki2022neural} underscores that the ultimate goal is the optimal dual solution, with all intermediate duals serving as approximations along the way. Consequently, when selecting from among various alternatives for $\widehat{u}_t^*$  at a given iteration $t$, our aim is to identify the one closest to $u^*$.  To accomplish this, in the formulation \eqref{spd:relint}-\eqref{spd:eq2}, we opt for the following alternative objective function instead of \eqref{spd:relint}.

\begin{align}
\min  \quad & |u^*-\widehat{u}_t^*| \label{spd:goal}
\end{align}

Notably, since we are dealing with a Mixed-Integer Programming (MIP) as our Master Problem (MP), $u^*$ may not be the sole contributor to convergence unless the linear hull represent a very tight approximation of the integer hull. Moreover, even if $u^*$ is known, an immediate proof of optimality is not guaranteed.

\subsection{Imitation Learning}
\label{subsection:imitlearning}

While considering all the extreme rays and extreme points in \eqref{mp:obj}-\eqref{mp:eq4} (particularly in \eqref{mp:eq2} and \eqref{mp:eq3}), it is noteworthy, as per \cite[Theorem 4]{magnanti1981accelerating}, that for the convex hull formulation of a mixed-integer program, a single Benders cut suffices to demonstrate optimality. This is only valid, however, if this single cut originates from the optimal subproblem dual solution. Yet, in cases where an integer hull formulation is not available, the iterative refinement of an outer approximation with Benders cuts will eventually lead to optimality.\\

In \cite{magnanti1981accelerating}, no systematic method for producing a core point is proposed. Conversely, in \cite{papadakos2008practical}, an \emph{approximate core point} is computed and employed in each iteration to generate non-dominated Benders cuts. The process is as follows: \cite{papadakos2008practical} initializes a core point approximation $x^0$ with the first solution to the MP. At iteration $t$, the subproblem dual is solved using the proposed solution from the master problem at that iteration. The approximation is then updated at each iteration $t$ by selecting a solution halfway between the previous  point (which is already cut off) and the solution of the MP for that iteration. While authors assert that  this gradual update of the core point estimation brings it closer to an actual core point, it remains unclear how connecting a point along the line from the integer hull (the new optimal) to the previous optimal solution(which is already cut off) can be belong (or converge) to a relative interior of integer hull of a given iteration.\\

\subsection{Procedure to Compute a Core Point}

Our approach operates as follows: We execute our Benders decomposition within a branch-and-bound framework. This entails activating the solution pool of our modern mixed-integer programming solver or record all visited integer solutions. Initially, we solve the master problem without incorporating the feasibility and optimality cuts:

\begin{align}
z^*_{MP} = \min_{x} \quad & c^T x + \eta \label{pool:obj}\\
\text{s.t.}  \quad & Fx \leq g \label{pool:eq1}\\
 \quad & \mathbf{x}  \in \mathbb{Z}\label{pool:eq4}
\end{align}

We enable the solver to populate the solution pool with a specified number of solutions for this problem, denoted as $N^{feas}$. If this quantity of solutions is insufficient, we employ an iterative approach. Specifically, we exclude the current optimal solution from the Master Problem (MP), which refers to the MP solved without the inclusion of feasibility and optimality cuts. Subsequently, we iteratively resolve this modified MP to attain a new optimal solution. This strategy involves a sequence of iterations, each excluding the solution found in the preceding step. In each of these iterations, the Master Problem is solved, and the resulting solution is recorded. Following this iterative process, a combinatorial cut $\sum_{j\in J^0} x_j + \sum_{j\in J^1}(1-x_j) \geq 1$  is appended to exclude the current solution, where $J^0$ and $J^1$ denote the sets of indices for variables taking values of 0 and 1, respectively, in the reported optimal solution. Subsequently, the master problem is further iteratively resolved to accumulate a sufficient number of integer solutions. All of these solutions reside within the integer hull of the initial MP. Leveraging Lemma \autoref{lem:relint}, a core point within the integer hull of \eqref{pool:obj}-\eqref{pool:eq4} can be ascertained. It's important to note that all the accumulated combinatorial cuts are to be removed prior to initiating the Benders process, once a sufficient number of solutions have been acquired.\\

In cases where the number of such solutions is adequately large, as Benders iterations commence, the solution pool can be periodically checked to record new solutions encountered throughout the process, possibly through the aid of MIP heuristics such as RINS, FeasOpt etc. Concurrently, previously recorded solutions that no longer pertain to the MP polytope at a given iteration can be filtered out. Given the existing optimal solution and at least one additional solution from the pool that remains feasible within the MP polytope of the present iteration, Lemma \autoref{lem:relint} provides another core point within the integer hull of solutions within the master problem for the current iteration.\\

\def\circledarrowclockwise#1#2#3{ 
\draw[#1,->] (#2) +(80:#3) arc(80:-260:#3);
}
\def\circledarrowcounterclockwise#1#2#3{ 
\draw[#1,->] (#2) +(10:#3) arc(10:350:#3);
}

%
%

This core point is subsequently utilized in \eqref{spd:relint}-\eqref{spd:eq2} to compute the duals. This leads to the compilation of state-action pairs as  $\mathcal{D}\{(s_i,a_i^*)\}_{i=1}^N$. These actions are presumed to represent the optimal actions executed by an expert, employing the aforementioned technique. Our goal is to learn such a policy.

\subsection{Policy Learning}

When calculating the loss function, a distinction arises between the two strategies outlined in \autoref{policy:MW}. In the first strategy, we compare our inferred duals with the observations for the same input. However, in the second strategy described in \autoref{policy:lid}, the comparison is limited to the optimal non-dominated duals, represented by $u^*$, at which convergence has been achieved.\\

To maximize the likelihood of the expert action, we minimize the following cross-entropy loss function:

\begin{align}
L(\eta) = -\sum_{p \in \mathcal{P}}\sum_{t=1}^{T^{p_i}} \log p_\theta (a_t^{*(i)}|s_t^{(i)}).
\end{align}
wherein $T^{p_i}$ represents the number of iterations taken to satisfy the termination criteria in instance  $p\in\mathcal{P}=\{1,2,\dots,  P\}$. This information will play a pivotal role in our policy learning process.
%
%


\subsection{State representation}

While there are alternative approaches, such as the bipartite graph representation of Mixed-Integer Programming (MIP) as discussed in \cite{Gasse2019}, we chose to adopt the methodology presented in \cite{ding2019accelerating}. We made this choice due to the alignment of the transition flow with our specific objectives. We represent the state with a tripartite graph denoted as  $\mathcal{G} = (\mathcal{V}^{MP}, \mathcal{C}^{MP}, \{o\}, \mathcal{E})$. In this graph, each node in $\mathcal{V}^{MP}$ corresponds to a variable in the master problem,  a node in $\mathcal{C}^{MP}$  corresponds to a constraint in the master problem, and an objective node $o$ is present. $\mathcal{E}$ comprises edges  $(i,j)$ following these rules: An edge $(v,c)$ appears when the corresponding variable is found in a constraint; an edge $(v,o)$  appears because all variables contribute in the objective function; and an edge $(c,o)$ is established as the right-hand side also contributes to the objective function. Of particular significance are the edge sets $(v,o)$ and $(c,o)$ as they directly link $\eta$ to the Benders cuts in  $\eqref{mp:eq2}$. In summary:

\begin{itemize}
  \item $\{c,v\}, \forall c\in \mathcal{C}^{MP}, v\in \mathcal{V}^{MP}$, if the corresponding variable has a non-zero coefficient in the respective constraint within the MP. The features include the non-zero entries in the coefficient matrix.
  \item $\{v,o\}, \forall v\in \mathcal{V}^{MP}$ exists between each variable and the objective function. Associated features encompass the variable's cost in the objective function.
  \item  $\{c,o\}, \forall v\in \mathcal{C}^{MP}$ is established between each constraint and the objective function node. As a result, the right-hand side of the constraint can be considered as a feature.
\end{itemize}

\subsection{Policy Parameterizations}

The policy for selecting dual values, denoted as $\pi(a|s_t)$ , can be effectively represented using a Graph Convolutional Network (GCN) architecture (\cite{Scarselli2009,Gori2005,Bruna2014}). In line with this representation, a graph convolutional layer is naturally proposed in \cite{ding2019accelerating} to process the tripartite graph. At each iteration $t$, the input to our model is a state $s_t =  \mathcal{G}$ represented as a graph, which undergoes the following transitions and information transformations: 1) $\mathcal{V}^{MP}$ transitions to  $\{o\}$. 2) $\mathcal{V}^{MP}\bigcup \{o\}$ transitions to $\mathcal{C}^{MP}$. 3) $\mathcal{C}^{MP}$ transitions to  $\{o\}$. 4)  $\mathcal{C}^{MP}\bigcup \{o\}$ transitions to  $\mathcal{V}^{MP}$. \\


Each node within this tripartite graph undergoes an embedding process, ensuring uniform embedding sizes across all node types. A single convolutional operation is divided into distinct steps for each instance of the graph $\mathcal{G}^{MP}$:

\begin{align}
&   o' \leftarrow f_{o} (o, \sum_{i: (i,o)\in \mathcal{E} } g_{o}( v_i, o, e_{i,o})) \label{gcn:conv:c}\\
&   c'_i \leftarrow f_{\mathcal{C}} (c_i, \sum_{j: (i,j)\in \mathcal{E} } g_{\mathcal{C}}(c_i, v_j, o', e_{i,j})+ g_{\mathcal{C}}(c_i, o', e_{o,i})) \label{gcn:conv:c}\\
&   o' \leftarrow f_{o} (o, \sum_{i: (i,o)\in \mathcal{E} } g_{o}( c'_i, o', e_{i,o})) \label{gcn:conv:c}\\
&   v'_j \leftarrow f_{\mathcal{V}} (v_j, \sum_{i: (i,j)\in \mathcal{E} } g_{\mathcal{V}}(c'_i, v_j, o', e_{i,j})+ g_{\mathcal{V}}(v_i, o', e_{o,i})) \label{gcn:conv:c}
\end{align}

wherein $f_{\mathcal{C}}, g_{\mathcal{C}},g_{o}, f_{\mathcal{V}}$, $g_{\mathcal{V}}$ and $g_{o}$ are perception layers. 

\subsection{Differentiable Optimization Layer}
Up to this stage, the outcome consists of a collection of embeddings for the Master Problem (MP) variables. Each embedding constitutes an array within the Rd $R^d$, where, in our context, $d=64$. Let $z^*_{SPD(u)}$ represent the optimal objective value for SPD, given the current input variable $\overline{x}$. We can construct an array $\bf q$ by extracting the first element from each embedding of the MP variables. This array is then utilized to minimize the negation of the objective function $u(b- A{ \bf q})$ over the SPD polytope. Here, $\bf u$ represents variables, and  $\bf q$ constitutes the learnable parameters. The output variable embeddings from the GCN layer should subsequently traverse through a fully connected layer and undergo ReLU activation. This process aims to propose a solution for the following problem:

\begin{align}
\min  \quad & -u(b- A{ q}) \\
\text{s.t.} \quad & u~B \geq d  \\
& z^*_{SPD(u)} =  u(b- A\overline{x}) \\
&   {u}\in \mathbb{R}^+\cup\{0\} .
\end{align}
for the \emph{Magnanti-Wong} policy. We can introduce $\bf P'=(P\overline{x})$ and $\bf Q=(P')^TP'$, where $\bf P$ represents the matrix of MP variable embeddings:

\begin{align}
\min  \quad & \frac{1}{2} u^T{\bf Q}u \\
\text{s.t.} \quad & u~B \geq d  \\
& z^*_{SPD(u)} =  u(b- A\overline{x}) \\
&   {u}\in \mathbb{R}^+\cup\{0\} .
\end{align}
for the \emph{Laster Iteration Duals (LID)} policy.\\

The differentiable optimization library OptNet \citep{amos2021optnet} offers functionalities that not only enable the solution of the aforementioned problem, but also automatically compute gradients with respect to the parameters.

\section{Computational Experiments}\label{sec::experiments}
The primary objective of our research is to assess the benefits of incorporating a learning process into the traditional Benders decomposition approach. Given that the OptNet solver operates solely on a single CPU, we are unable to leverage extensive parallelization. Consequently, we rely on a single CPU for solving each instance.\\

Our benchmark comprises a collection of instances of the Uncapacitated Single Allocation $p$ Hub Location Routing Problem (SAHLRP), exhibiting certain similarities to the work in \cite{DANACH2019597}, as detailed in \autoref{appendix}. This benchmark encompasses instances with 30, 40, 50, 60, and 70 nodes. In total, 1500 instances are randomly generated, with each instance size having 500 instances. Among these instances, 250 are allocated for training, 150 for validation, and the remaining 100 for testing. The models are trained using the Adam optimizer, employing a learning rate of $3\times 10^{-4}$, and this learning rate is divided by 5 if no improvement is observed in the validation error reduction. Mini-batches of size 32 are employed during training. When executing the regular Benders decomposition, a single thread is employed, and the computation is performed on an Intel i9 processor with a clock speed of 2500 MHz and 32 GB of RAM.\\

The open-source framework \texttt{BENMIP--Benders for Mixed Integer Programming} has been developed and maintained over the course of several years. This framework is written in C++, with a provided Python interface, and can be compiled for both Windows (using the Microsoft C++ compiler) and Linux (specifically Ubuntu) platforms. For the current project, we compiled \texttt{BENMIP} using the VC++ (Visual Studio 2019) compiler on a system running Ubuntu. The MIP solver utilized is CPLEX 20.1.0.0.

Both the Master Problem (MP) and the Subpoblem Dual (SDP) are solved using the CPLEX solver. Specific CPLEX parameters have been selected and configured for the solving process. For instance, \texttt{IloCplex::MIPEmphasis} is set to \texttt{CPX\_MIPEMPHASIS\_BALANCED}. Custom termination criteria for CPLEX, specifically \texttt{UserAbort}, are employed under the following conditions: CPLEX terminates if a maximum of 5 hours (equivalent to 18000 seconds) of CPU time has been utilized and an incumbent solution has been identified. In terms of computational efficiency, $N^{feas}$ is set to a minimum of 10,000 instances. Additionally, a rigorous process of extensive preliminary computational experiments has led to the decision to maintain an optimality tolerance of 0.5 percent.

\subsection{Training Data Generation}

The approach described in \autoref{subsection:imitlearning} is employed to ensure that a point consistently exists within the closure of the relative interior of the MP integer hull. All instances are solved to optimality (given the tolerance), and the recorded data includes the state of the master problem at each iteration, the optimal solution of the MP at iteration $t$ denoted as $(\overline{x}_t,\overline{\eta}_t)$, the relative interior point of that specific iteration, and the dual solution obtained using the Magnanti-Wong method. Specifically, for the second approach, the dual solutions selected are those upon which the overall optimality has been proven.

\subsection{Analysis}

For each individual instance, we utilize the \emph{Magnanti-Wong} method to solve the problem. During this process, we record various indicators, including the computational time required and the number of calls to the separation routine, or equivalently, the number of iterations taken to meet the termination criteria. It is worth mentioning that even if the subproblem can be decomposed into smaller problems, we choose not to do so. Instead, we solve it as a single problem. This decision is made to ensure a fair comparison with the neural-based techniques and serves as the baseline for our study.\\

While it is accurate to state that the original model consists of a substantial number of variables and constraints, it is important to note that the growth rate of these numbers remains limited, specifically on the order of $\mathcal{O}(n^2)$. Additionally, the sparsity of the graph and the manageability of the GCN play crucial roles in keeping the computations within a feasible range, even as more cuts are incorporated into the process.\\

\begin{figure}[htbp]
    \centering
    \includegraphics[width=0.9\columnwidth]{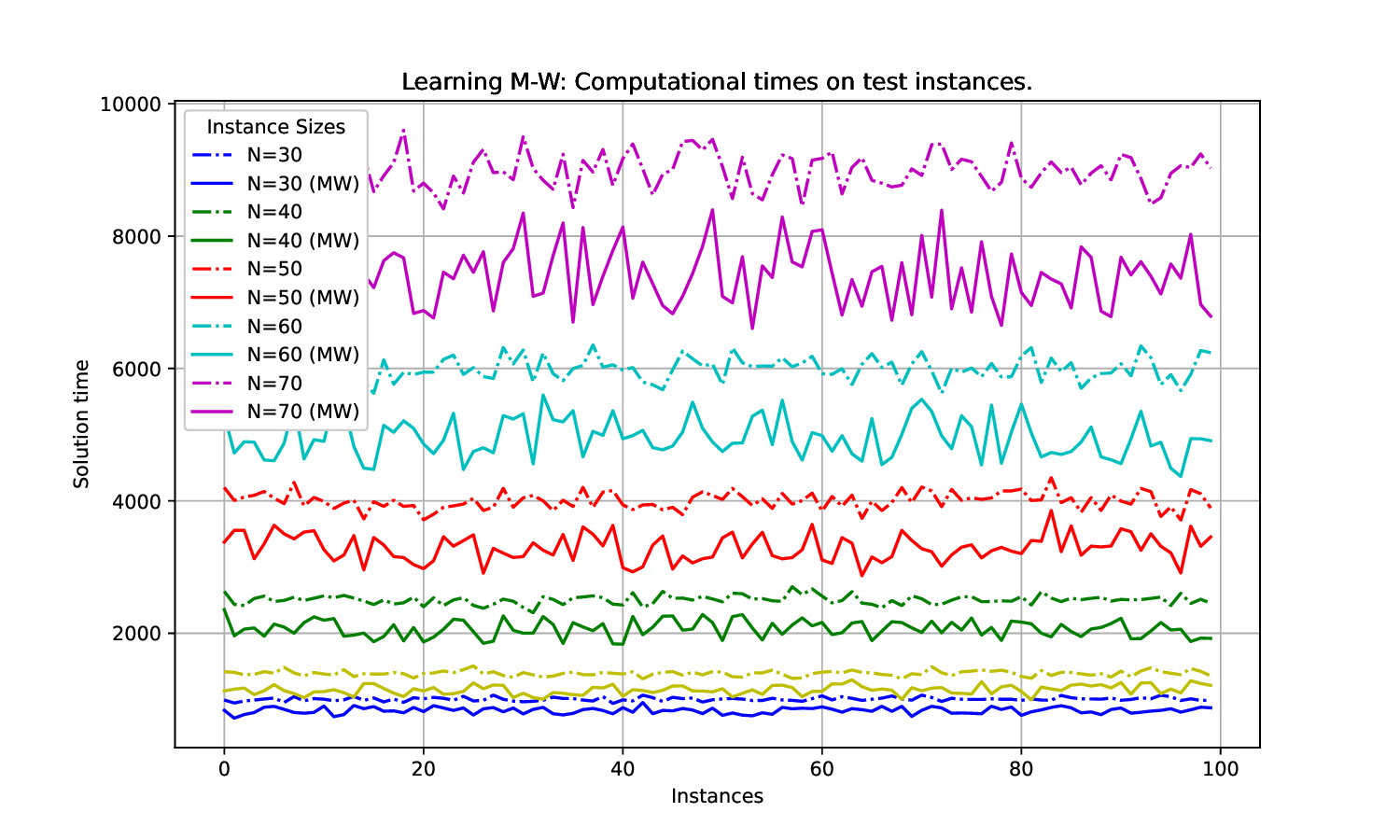}  
    \caption{BD-MW vs. M-W: Computational results for instances of different size.}
    \label{fig:ComputationalTimeMW}
\end{figure}

With regards to the results, we are now reporting the outcomes based on 100 test instances. In \autoref{fig:ComputationalTimeMW}, we commence by comparing the computational times required for the M-W Benders decomposition with those for the learned M-W policy. It is evident from the graph that the learned M-W policy consistently establishes a lower envelope in comparison to the M-W Benders decomposition method.\\

\begin{figure}[htbp]
    \centering
    \includegraphics[width=0.9\columnwidth]{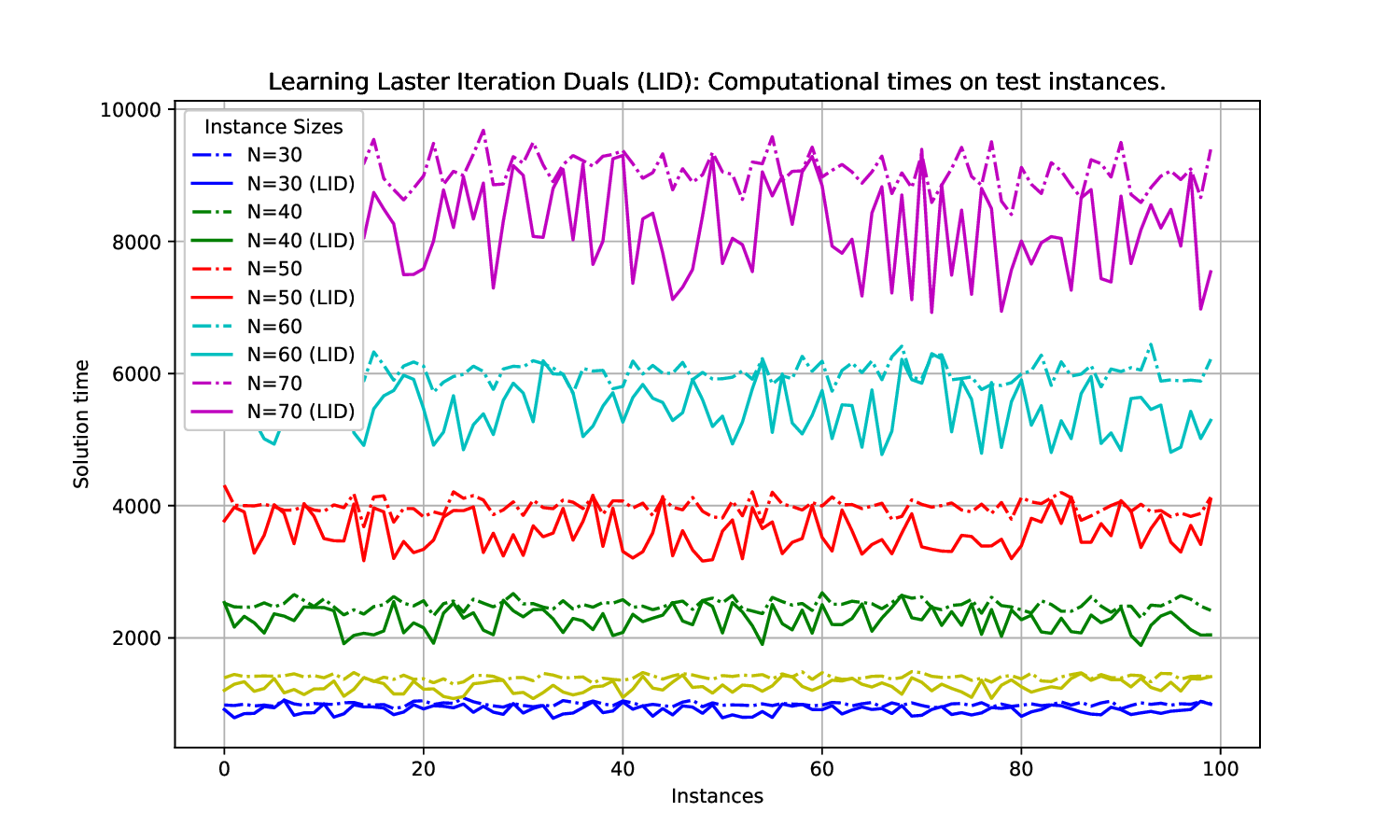}  
    \caption{BD-MW vs. LID: Computational results for instances of different size.}
    \label{fig:ComputationalTimeLID}
\end{figure}

However, in contrast to the clear lower envelope seen in the M-W approach, the lower bound is not as apparent in the LID approach, as depicted in \autoref{fig:ComputationalTimeLID}. While it is evident that the LID learned policy frequently outperforms the BD-MW approach, this absolute superiority isn't consistently established.\\

\begin{figure}[htbp]
    \centering
    \includegraphics[width=0.9\columnwidth]{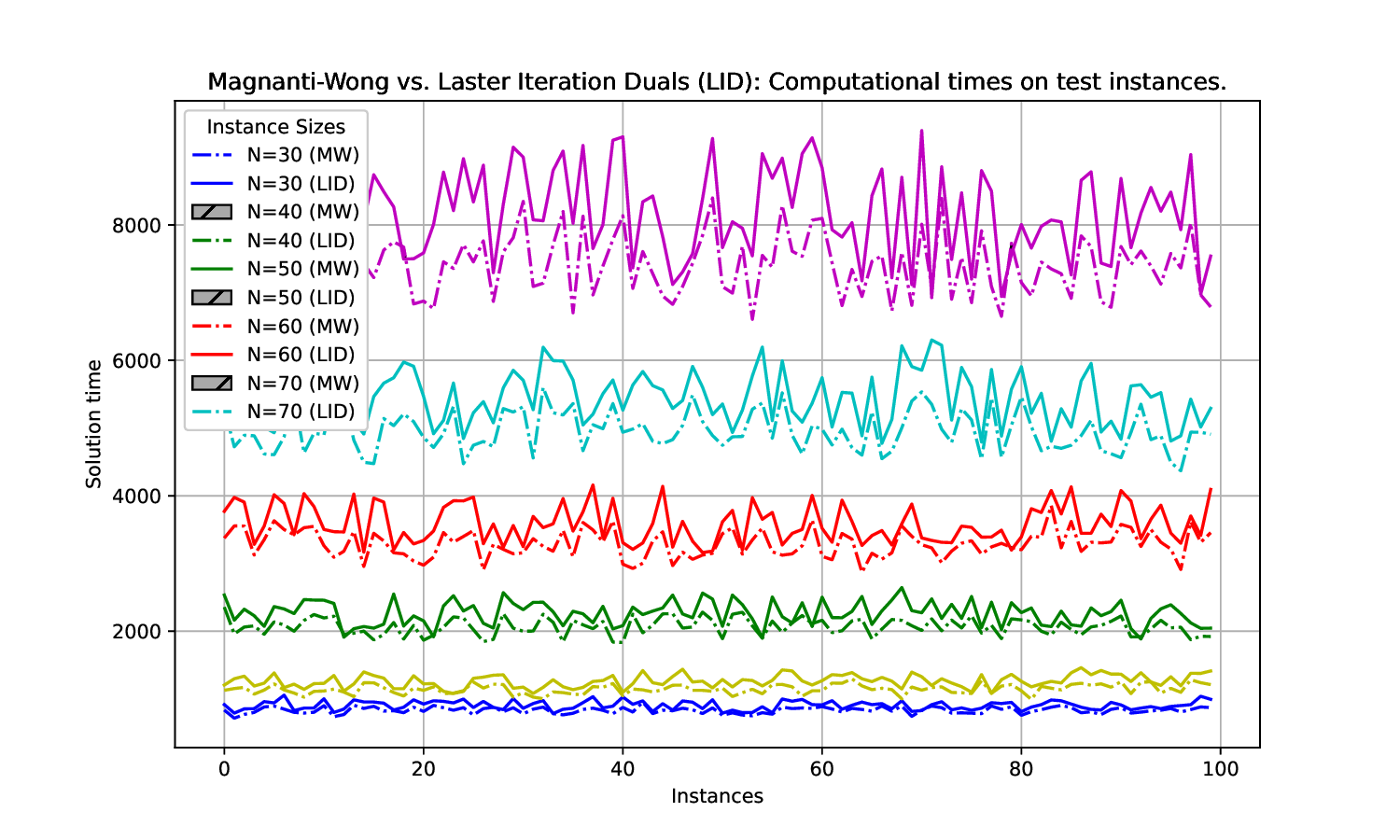}  
    \caption{MW vs. LID: Computational results for instances of different size.}
    \label{fig:ComputationalTimeMWLID}
\end{figure}

It is evident from the comparison presented in \autoref{fig:ComputationalTimeMWLID} that the learned M-W policy consistently outperforms the learned LID policy.

\section{Summary, Conclusion and an outlook to future research directions}\label{sec::conclusion}

In this study, we introduced two novel approaches to enhance the convergence of Benders decomposition, particularly in scenarios where the subproblem exhibits primal degeneracy and multiple dual optima. Our computational experiments yielded results that affirm the efficacy of the learned Magnanti-Wong policy, surpassing the performance of the traditional Magnanti-Wong technique in Benders decomposition (BD-MW). While the learned Last Iteration Dual policy also exhibited improvements over BD-MW, it did not consistently establish a lower envelope in computational time. The overall performance of the proposed approaches is highly satisfactory, with more than 80\% of instances showcasing superiority. An interesting aspect is that, with learned policies, the need for explicit introduction of the relative interior is obviated, streamlining the complex aspect of the process.\\

Future investigations should delve into the inclusion of feasibility cuts in the proposed frameworks. Moreover, exploring alternative methods for solving the quadratic model is crucial to alleviate the computational bottleneck inherent in the use of a single CPU in OptNet.


%

\section*{Acknowledgement}
Sharkey Predictim Globe\footnote{\href{www.predictim-globe.com}{www.predictim-globe.com}} company has provided us with the necessary resources to conduct our computational experiments and development.  Additionally, the BENMIP project, funded by PGMO - the Gaspard Monge Program for Optimization and Operational Research, has been instrumental in our research. BENMIP framework has been extended to incorporate the new features.

\bibliographystyle{elsarticle-harv}
\bibliography{introBib1,petroleum,crossdock,MIP,vns,alns,BendersDecomposition,gnn,deeplearning,ai}

\clearpage
\appendix
\label{appendix}

\section{Uncapacitated Single Allocation $p$ Hub Location Routing Problem (USA$p$HLRP)}

The benchmark used in our study consists of instances of the \emph{Uncapacitated Single Allocation $p$ Hub Location Routing Problem (USA$p$HLRP)}. This problem is derived from the original model introduced in \cite{DANACH2019597} for the \emph{Capacitated Single Allocation $p$ Hub Location Routing Problem (CSA$p$HLRP)}. In our problem, we have eliminated the capacity constraints and other design constraints to create a simplified version. The relevant variables and parameters used to formulate this problem are presented in  \autoref{tab:Params} and \autoref{tab:DecVars}.\\
\begin{table}[!h]
\begin{center}
 \caption{Model Parameters.}
\begin{tabular}{p{0.15\columnwidth}p{0.7\columnwidth}}
\toprule
$w_{ij}$: &The flow from $i$ to $j$,\\
$t_{ij}$: &The distance or time on a direct link (arc) between nodes  $i - j$,\\
$\alpha$: &The factor representing economies of scale, indicating the efficiency of travel time over hub edges, \\
$p$: &The upper bound on the number of hubs or depots,\\
$q$: &The lower bound on the number of hubs or depots,\\
$\Gamma$: &The minimum number of spokes allocated to each hub or depot node,\\
$\varphi^k$: &The fixed average transshipment time at hub or depot node $k$,\\
$I_{kl}$: &The fixed cost of allocating spoke node $k$ to the hub node $l$,\\
$F_k$: &The fixed cost of upgrading node  $k$ to a hub node.
 \\\bottomrule
\end{tabular}   \label{tab:Params}
\end{center}
\end{table}

\begin{table}[!h]
 \begin{center}
 \caption{Decision Variables.}
\begin{tabular}{p{0.15\columnwidth}p{0.7\columnwidth}}
\toprule
$x_{ijkl}$:  &The fraction of flow from $i$ to $j$ traversing inter-hub edge $\{k,l\}$,\\
$s_{ijkl}$:  &The fraction of flow from $i$ to $j$ traversing non-hub edge $\{k,l\}$,\\
$z_{ik}$: &1,  if node $i$ is allocated to node $k$ where $k$ is a hub, 0 otherwise.
\\\bottomrule
\end{tabular}   \label{tab:DecVars}
\end{center}
\end{table}

(SAHLRP)

\begin{align}
\min \:      &
            \sum_{i,j,k,l} \left(t_{kl}(s_{ijkl}+{\alpha x_{ijkl}})\right) + \nonumber\\ 
            & \sum_{i,j,k,l:(k\neq i \vee l\neq j)  }(\varphi^k +  \varphi^l)x_{ijkl} +\nonumber\\
            & \sum_k F_k z_{kk} + \sum_{k,l} I_{kl} z_{kl}\label{eq:obj}
            \\
s.\:t.\:\nonumber\\
            &   \sum_l z_{kl}  =1         &       \forall k \in V   \label{eq:eq2}\\
            &    z_{ik} \leq z_{kk}         &       \forall i,k\in V: k\neq i    \label{eq:eq4}\\
            &   \sum_{i} z_{ik} \geq \Gamma z_{kk}          &       \forall k \in V   \label{eq:eq5}\\
             & \sum_{k \neq i} (x_{ijik} + s_{ijik}) = 1, &         \forall i,j\in V:j\neq i,
        \label{eq:eq19}\\
            &   \sum_{l \neq j} (x_{ijlj} + s_{ijlj}) = 1, &       \forall i,j\in V:j\neq i,
        \label{eq:eq20}\\
            &\sum_{l \neq i,k} (x_{ijkl} + s_{ijkl}) = \nonumber\\
            & \sum_{l \neq j,k} (x_{ijlk} + s_{ijlk}), &         \forall i,j,k  \in V,\: k \not \in \{ i,j,\}\:,
        \label{eq:eq21}\\
            &   \sum_{l\neq k} x_{ijkl} \leq z_{kk}    &       \forall i,j,k \in V:j\neq i, k<l \label{eq:eq22}\\
             &   \sum_{l\neq k} x_{ijlk} \leq z_{kk}    &       \forall i,j,k \in V:j\neq i, k<l \label{eq:eq23}\\
		    %
            & z \in \mathbb{B}^{|V|\times |V|}, ~~x_{ijkl}, s_{ijkl} \in \in \mathbb{R}^{|V|^4}_{[0,1]} 
\end{align}

The objective function \eqref{eq:obj} consists of four parts: the first part represents the total transportation time, the second part accounts for transshipment time, the third part consists of fixed costs for upgrading a node to a hub node, and the final part deals with the assignment of nodes to hub nodes. The transportation component considers travel time on both hub-level edges and spoke-level arcs. The travel time on hub-level edges is discounted by the factor $0 < \alpha < 1$, reflecting the faster services offered by transporters. Transshipment time is taken into account for hub nodes along each O-D path.\\

Constraints \eqref{eq:eq2} ensure that every node is assigned to one and only one hub. A self-allocation at node $i$ indicates that node $i$ is a hub (i.e., $z_{ii}=1$). Constraints \eqref{eq:eq4} ensure that a spoke node can only be assigned to a hub node. Additionally, constraints \eqref{eq:eq5} guarantee that at least $\Gamma$ nodes are allocated to each designated hub node.

Constraints \eqref{eq:eq19} - \eqref{eq:eq21} correspond to flow conservation for every origin-destination (O-D) pair. Constraints \eqref{eq:eq22} and \eqref{eq:eq23} ensure that traversing an edge composed of two hub nodes is equivalent to traversing a hub edge connecting those two hubs.


The proposed Benders decomposition for this problem assumes that all integer variables remain in the master, including the constraints \eqref{eq:eq2}-\eqref{eq:eq5}. The continuous variables, accompanied by constraints \eqref{eq:eq19}-\eqref{eq:eq23}, constitute the subproblem.

\end{document}